\documentstyle[theorem]{article}

\theoremstyle{changebreak}
\theorembodyfont{\rm}
\newtheorem{sdef}{Definition}[section]
\newtheorem{sprop}[sdef]{Proposition}
\newtheorem{sthm}[sdef]{Theorem}
\newtheorem{slem}[sdef]{Lemma}
\newtheorem{scor}[sdef]{Corollary}
\newtheorem{sques}[sdef]{Question}
\newcommand{\Db}[2]{{{\bf \Delta}^#1_#2}}
\newcommand{\Sb}[2]{{{\bf \Sigma}^#1_#2}}
\newcommand{\Dl}[2]{{{\Delta}^#1_#2}}
\newcommand{\Sl}[2]{{{\Sigma}^#1_#2}}
\newfont{\Euler}{msbm10}
\newfont{\Fraktur}{eufm10}
\newcommand{\HechL}[1]{{{\sf Hech}({\bf L}[#1])}}
\newcommand{\CohL}[1]{{{\sf Coh}({\bf L}[#1])}}
\newcommand{\RanL}[1]{{{\sf Ran}({\bf L}[#1])}}
\newcommand{\BBB}{{\mbox{{\Euler B}}}}
\newcommand{\CCC}{{\mbox{{\Euler C}}}}
\newcommand{\DDD}{{\mbox{{\Euler D}}}}
\newcommand{\NC}[1]{{{\sf N}_{\CCC}({\bf L}[#1])}}
\newcommand{\ND}[1]{{{\sf N}_{\DDD}({\bf L}[#1])}}
\newcommand{\BCC}[1]{{{\sf BC}_{\CCC}({\bf L}[#1])}}
\newcommand{\BCD}[1]{{{\sf BC}_{\DDD}({\bf L}[#1])}}
\newcommand{\LLL}{{\mbox{{\Euler L}}}}
\newcommand{\MMM}{{\mbox{{\Euler M}}}}
\newcommand{\PPP}{{\mbox{{\Euler P}}}}

\newcommand{\SSS}{{\mbox{{\Euler S}}}}

\newcommand{\SiiD}{\Sb{1}{2}(\DDD)}
\newcommand{\SiiC}{\Sb{1}{2}(\CCC)}
\newcommand{\DiiD}{\Db{1}{2}(\DDD)}
\newcommand{\DiiC}{\Db{1}{2}(\CCC)}
\newcommand{\SiiL}{\Sb{1}{2}(\LLL)}
\newcommand{\SiiM}{\Sb{1}{2}(\MMM)}
\newcommand{\DiiL}{\Db{1}{2}(\LLL)}
\newcommand{\DiiM}{\Db{1}{2}(\MMM)}
\newcommand{\SiiS}{\Sb{1}{2}(\SSS)}
\newcommand{\DiiS}{\Db{1}{2}(\SSS)}
\newcommand{\frM}{{\mbox{{\Fraktur M}}}}
\newcommand{\frnM}{{\mbox{{\Fraktur M}}}}
\newcommand{\omlom}{\omega^{<\omega}}
\newcommand{\omom}{\omega^\omega}

\newcommand{\stem}{{\mbox{stem}}}
\newcommand{\Split}{{\mbox{Split}}}
\newcommand{\Succ}{{\mbox{Succ}}}
\newcommand{\mod}{{\mbox{mod}}}
\newcommand{\code}{{\mbox{{\sf code}}}}

\newcommand{\et}{\,{}\hat{}\,}
\newcommand{\la}{\langle}
\newcommand{\ra}{\rangle}
\newcommand{\ran}{{\mbox{ran}}}
\newcommand{\dom}{{\mbox{dom}}}

\newcommand{\qed}{\begin{flushright}{\sf q.e.d.}\end{flushright}}

\begin{document}

\begin{center}

{\LARGE Solovay--type characterizations\\ 
for forcing--algebras}\\[0.5cm]
{\large J\"org Brendle}\footnote{Part of this research was done
while the first author was supported by DFG--grant Nr. Br 1420/1--1
and the second author by DAAD--grant Ref.316--D/96/20969
in the program HSP II/AUFE and a grant of the Studienstiftung
des Deutschen Volkes.\\[0,5cm]
}\\
Department of Mathematics\\
Dartmouth College, Hanover, NH 03755, USA\\
{\footnotesize Brendle@MAC.dartmouth.edu}\\[0.5cm]
{\large Benedikt L\"owe}\\
Department of Mathematics\\
University of California, Berkeley, CA 94720, USA\\
{\footnotesize loewe@math.berkeley.edu}\\
\today\\
\end{center}

\begin{abstract}
We give characterizations for the (in ZFC unprovable) sentences
``Every $\Sb{1}{2}$--set is measurable" and 
``Every $\Db{1}{2}$--set is measurable" for various notions
of measurability derived from well--known forcing partial orderings.
\end{abstract}

%%%%%%%%%%%%%%%%%%%%   Introduction/Section 1   %%%%%%%%%%%%%%%%%%%

\section{Introduction}
In recent years, forcing notions which were originally devised to carry
out some consistency proof have emerged more and more as independent
mathematical objects which should be studied in their own right,
from various angles.
One such endeavor has been to investigate notions of measurability 
(that is, $\sigma$--algebras) associated 
with  forcing orderings adding a generic real.
This has a long tradition since the notions related to
{\sc Cohen} and random forcing are  the {\sc Baire} property
and {\sc Lebesgue} measurability which have always been
in the focus of set--theoretic research
({\it cf.} the results of
\cite{S70} and \cite{Sh84}). Other algebras which have
been around for quite a while include the {\sc Marczewski}--measurable
sets  \cite{Mar} which correspond to {\sc Sacks} forcing and
the completely {\sc Ramsey} sets which are connected with 
{\sc Mathias} forcing. In all of these cases,   measurability
of the analytic sets has been proved long ago,
and it has been known that one can get non--measurable sets on
the $\Db{1}{2}$--level in the constructible universe
{\bf L}. Furthermore, {\sc Solovay} (see \ref{SoloSb12})
proved in the sixties that the statement ``all $\Sb{1}{2}$--sets
are {\sc Lebesgue}--measurable" is equivalent to ``over each ${\bf
L}[a]$, there is a measure--one set of random reals" which is in turn
equivalent to ``for all $a$, the union of all null sets coded in
${\bf L}[a]$ is null", thereby reducing a statement about measurability of
projective sets to what might be termed a {\it transcendence principle}
over ${\bf L}$. The value of such characterizations, apart from their
intrinsic beauty, is obvious: they make it much easier to check
whether $\Sb{1}{2}$--measurability holds in a given model of set theory.
So it is a natural question whether statements like
``all $\Db{1}{2}$--sets are 
$\PPP$--measurable" and ``all $\Sb{1}{2}$--sets are 
$\PPP$--measurable" can be characterized in {\sc Solovay}'s
fashion as  transcendence
principles over ${\bf L}$, for other forcing notions $\PPP$ adding 
a generic real.

In this work, we show this can be done in several cases. 
The most interesting results concern {\sc Hechler} forcing $\DDD$,
the standard c.c.c. forcing notion adjoining a dominating real,
and the related {\it dominating topology} ${\cal D}$ on $\omom$
(see the definition in \ref{Maindefi}, (ii)).
The notion of measurability associated with $\DDD$ is,
of course, the property of {\sc Baire} with respect to ${\cal D}$.
We show that  all $\Db{1}{2}$--sets have the {\sc Baire} property
in ${\cal D}$ iff all $\Sb{1}{2}$--sets have the {\sc Baire}
property in the standard topology on $\omom$ (Theorem \ref{DiiD}).
Using a combinatorial result on the dominating topology
due to \cite{LR95} which builds, in turn, on the
combinatorics of {\sc Hechler} forcing developed in \cite{BJS92},
we then get, as a rather easy consequence of the characterization
on the $\Db{1}{2}$--level, that all $\Sb{1}{2}$--sets have the {\sc Baire} 
property in ${\cal D}$ iff $\aleph_1^{{\bf L}[a]} <
\aleph_1$ for all reals $a$ (Theorem \ref{HechSigma}).
This confirms a conjecture put forward by {\sc Judah}
(private communication). It's the only case we know of where
the consistency strength of $\Sb{1}{2} - \PPP$--measurability
is already an inaccessible. This should be compared to 
the result of \cite{Sh84} showing that the consistency
strength of $\Sb{1}{3}$--{\sc Lebesgue}--measurability is
an inaccessible.

We also investigate various other notions of measurability,
e.g. $\MMM$--measura\-bi\-lity which is derived from {\sc Miller}'s
rational perfect set forcing $\MMM$. We show that all $\Db{1}{2}$--sets
are $\MMM$--measurable iff all $\Sb{1}{2}$--sets are iff
$\omom\cap {\bf L}[a]$ is not dominating in $\omom$ for
all reals $a$ (Theorem \ref{MCharac}). In all the cases we
consider here, the proof of the projective statement assuming
the transcendence principle follows either from known
game--theoretic arguments or by rewriting the corresponding
proof for the standard {\sc Baire} property. Our main technical
results (\ref{Laverchar},   \ref{MCharac}, but also
\ref{DiiDbeschr} which follows from \ref{Laverchar} and \ref{DiiDDiiL}), then,
deal with the other direction --- taken care of by a {\sc Fubini}--argument
in case of the {\sc Baire} property and {\sc Lebesgue}
measurability which does not apply in our case ---
and have all a similar flavour: each time, we construct
a $\Dl{1}{2} (a)$--partition of the reals along a carefully
chosen scale of ${\bf L}[a]$. For example, to prove 
Theorem \ref{MCharac} mentioned above, we produce, under the
assumption that $\omom \cap {\bf L} [a]$ is dominating,
a $\Dl{1}{2} (a)$--super--{\sc Bernstein} set, where 
$A \subseteq \omom$ is called super--{\sc Bernstein}
iff both $A$ and $\omom\setminus A$ meet every superperfect set.

This paper is organized as follows. In section \ref{Maindefi}, we introduce
the notions of forcing we are interested in, define what we
mean by the corresponding notion of measurability and fix
our notation. Section \ref{General} contains general results 
on the connection between the various measurability notions we study. The next
three sections contain the main results: in section \ref{Hechsec}  we study
{\sc Hechler} forcing; sections 
\ref{Laversec} and \ref{Millersec} deal with {\sc Laver}
and {\sc Miller} forcing, respectively. We conclude with a brief
remark about {\sc Sacks} forcing in section \ref{Sackssec}, and
an overview on our results as well as an open problem in section \ref{Summ}.
All sections depend on sections \ref{Maindefi} and \ref{General}; the
higher--numbered sections can
be read independently of each other; however, \ref{DiiDbeschr}
uses  \ref{Laverchar}.

%%%%%%%%%%%%%%%%%%%%%%%%%%%   Section 2   %%%%%%%%%%%%%%%%%%%%%%%%%%%

\section{Main definitions and notation} \label{Maindefi}

\begin{enumerate}
\item $\CCC:=\la\omlom,\supseteq\ra$ is called {\sc Cohen} forcing.
For each condition $s$ we define $[s] :=\{f\in\omom:s\subseteq f\}$.
The sets $([s])_{s\in\omlom}$ are a topology base of the so called
{\sc Baire} space whose topology we denote by ${\cal B}$. Sometimes it
may be necessary to regard {\sc Cohen} forcing on the {\sc Cantor}
space $2^{\omega}$. In this case we will denote the topology by ${\cal C}$.
\item We call $\DDD:=\omega\times\omom$ {\sc Hechler} forcing,
when we have the following partial ordering on it:
\[\la N,f\ra\leq\la M,g\ra\iff N\geq M, f|M=g|M, f\geq g\]
We put $[N,f]:=\{x\in\omom:\; f|N \subseteq x$ and $x(n) \geq f(n)$
for all $n\}$.
Again, the sets $([N,f])_{\la N,f\ra\in\DDD}$ are a topology base of the 
{\it dominating topology} ${\cal D}$.
Obviously the dominating topology is finer than the {\sc Baire}
topology, because if we define the following real number
\[x_s(n):=\left\{\begin{array}{rl}s(n)&\mbox{ if } n<|s|\\
0&\mbox{else}\end{array}\right.\]
then $[|s|,x_s]=[s]$.
As we know from \cite{LR95}, ${\cal D}$ is a c.c.c. {\sc Baire} space.
\end{enumerate}
In contrast to these two forcings whose conditions form a topology
base on $\omom$ (and which we call therefore {\it topological
forcings}) we consider the following three {\it non--topological forcings}:
\begin{enumerate}
\setcounter{enumi}{2}
\item A tree $L\subseteq\omlom$ is called {\sc Laver} tree, if all
nodes above the stem are $\omega$--splitting 
nodes\footnote{A node is called {\it splitting} if it has more
than one immediate successor, and it is called $\omega$--{\it splitting} if
it has infinitely many immediate successors}.
We call the set of all {\sc Laver} trees ordered by inclusion
{\sc Laver} forcing $\LLL$.
\item A tree 
$M\subseteq\omlom$ is called {\it superperfect}, if
every splitting node is an $\omega$--splitting node and every node
has a (not necessarily immediate) successor which is a splitting
node (and therefore an $\omega$--splitting node).
{\sc Miller} forcing $\MMM$ is the set of all superperfect trees
ordered by inclusion.
\item In analogy to the definition of $\MMM$ we call a tree
$P\subseteq 2^{<\omega}$ {\it perfect}, 
if below every node there is a splitting node and define {\sc Sacks}
forcing $\SSS$ to be the set of all perfect trees ordered by inclusion.
\end{enumerate}
Given a tree $T \subseteq\omlom$, let $[T] := \{ f\in\omom :
\; f | n \in T$ for all $n\in\omega\}$ denote the set of its
branches. For $s\in T$, let Succ($s$) be the set of immediate
successors of $s$ in $T$. Split($T$) stands for the set of
splitting nodes of $T$.

We can associate each of these forcing in a natural
way with a notion of measurability.
In the definition of the topological forcings $\CCC$ and $\DDD$,
we remarked that the forcings form topology bases for ${\cal B}$ and
${\cal D}$ respectively. The forcings are therefore
quite naturally connected to the $\sigma$--algebra of sets with the
{\sc Baire} property in these topologies. The ${\cal B}$-- and 
${\cal D}$--meager sets are also called $\CCC$-- and $\DDD$--null sets.

In the case of non--topological forcings $\PPP\in\{\SSS,\MMM,\LLL\}$ we 
define a set of
real numbers $A$ ($A\subseteq\omom$ or $A\subseteq 2^\omega$ according to the 
definition of $\PPP$) to be {\it $\PPP$--measurable} if 
\[\forall p\in\PPP\;\exists p'\leq p\;([p']\cap A=\emptyset\mbox{ or }
[p']\cap \omom\setminus A=\emptyset)\]
and to be {\it $\PPP$--null} if
\[\forall p\in\PPP\;\exists p'\leq p\;([p']\cap A=\emptyset)\]
The ideal of $\PPP$--null sets we denote by $(p^0)$ and the set of 
complements of $\PPP$--null sets we denote by $(p^1)$.

For pointclasses $\Gamma$, we
abbreviate the sentence ``every set in $\Gamma$  is $\PPP$--measurable"
by $\Gamma(\PPP)$.
In addition to that we define a set $A$ to be 
{\it weakly $\PPP$--measurable} if either $A$ or its complement contains
the branches through some element of $\PPP$. As above, 
we abbreviate the sentence ``every
set in $\Gamma$ is weakly $\PPP$--measurable"
by $w\Gamma(\PPP)$. We will call a pointclass $\Gamma$ 
{\it topologically reasonable} if it is closed under continuous
preimages and has the following property:
$$\mbox{For }A\in\Gamma\mbox{ and }Q\mbox{ closed, we have }
A\cap Q\in\Gamma$$
\begin{slem}\label{weakM}
Let $\PPP$ be any of the forcings considered in this
work, and let $\Gamma$ be a
topologically reasonable pointclass. Then the following are equivalent:
\begin{enumerate}
\item $w\Gamma(\PPP)$
\item $\Gamma(\PPP)$
\end{enumerate}
\end{slem}
{\bf Proof :}\\
As the backward direction is obvious, we prove the forward direction:
Suppose $\Gamma(\PPP)$ is false. Then there is an $A\in\Gamma$ which is
not $\PPP$--measurable, {\it i.e.} there is a $P\in\PPP$ such
that for all $Q\leq P$:
$$[Q]\cap A\neq\emptyset$$
and
$$[Q]\cap \omom\setminus A\neq\emptyset$$
Let $\sigma$ be an homeomorphism between $[P]$ and $\omom$ (or
$2^\omega$ in the case that $\PPP$ is defined on the {\sc Cantor} space).
Then because of the properties postulated for $\Gamma$, $A\cap [P]$
and $A':=\sigma(A\cap [P])$ are in $\Gamma$. Because of $w\Gamma(\PPP)$
we have $Q'\in\PPP$ with either $[Q']\subseteq A'$ or $[Q']\subseteq 
\omom\setminus
A'$. Applying $\sigma^{-1}$ yields:
$$\sigma^{-1}[Q']\subseteq A\cap[P]\mbox{ or }
\sigma^{-1}[Q']\subseteq [P]\cap(\omom\setminus A)$$
But this is a contradiction.\qed
For this lemma we do not need closure under continuous
preimages, closure under homeomorphism would suffice.\\[0.5cm]
Let $\frnM$ be a model of ZFC and {\sf BC} be a fixed coding of the
{\sc Borel} sets. For a code $c$ we denote the decoded set with $A_c$.
If $\PPP$ 
is any of the defined forcing notions then
${\sf BC}_\PPP(\frnM)$ denotes the set of all real numbers in $\frnM$ which
code a Borel $\PPP$--null set (we have enough absoluteness properties
for this to make sense). For abbreviation we define:
$${\sf N}_\PPP(\frnM):=\bigcup\{A_c : {c\in{\sf BC}_\PPP(\frnM)}\}$$
For the c.c.c. forcings considered here,
one can prove that the $\PPP$--generic reals over
$\frnM$ are exactly those not in ${\sf N}_\PPP(\frnM)$.
This result allows the following definition:
\begin{sdef}
\begin{enumerate}
\item $\mbox{\sf{Coh}}(\frnM):=\omom\setminus{\sf N}_\CCC(\frnM)$
\item $\mbox{\sf{Hech}}(\frnM):=\omom\setminus{\sf N}_\DDD(\frnM)$
\end{enumerate}
\end{sdef}
\begin{sdef}
Let $A\subseteq\omom$. A real $x\in\omom$ is called
\begin{itemize}
\item {\it unbounded} over $A$, if:
$$\forall a\in A\;\exists^\infty n: a(n)\leq x(n)$$
\item {\it dominating} over $A$, if:
$$\forall a\in A\;\forall^\infty n: a(n)\leq x(n)$$
\end{itemize}
Shortly we write
$x\leq^*y:\iff \forall^\infty n:x(n)\leq y(n)$ for
``$y$ dominates $x$".\\
Let $B\subseteq\omom$. $B$ is called
\begin{itemize}
\item $\sigma$--{\it bounded} in $A$, if there is an $a\in A$
dominating over $B$
\item {\it unbounded} in $A$, if it is not
$\sigma$--bounded in $A$
\item {\it dominating} in $A$, if no $a \in A$ is unbounded over $B$
\end{itemize}
\end{sdef}
Apart from this, we use standard notions and notation of Descriptive
Set Theory and Forcing Theory (see e.g. \cite{Jech} or \cite{BJ95}).

%%%%%%%%%%%%%%%%%%%%   section 3   %%%%%%%%%%%%%%%%%%%%%%%%%%%%%%%%%%%

\section{General Results} \label{General}

We provide a few results on the connection between some
notions of measurability which hold for arbitrary topologically
reasonable pointclasses $\Gamma$.
\begin{sthm}\label{DiiDDiiC}
For any topologically reasonable pointclass $\Gamma$,
$\Gamma(\DDD)$ implies $\Gamma(\CCC)$
\end{sthm}
{\bf Proof :}\\
We define a mapping $\varphi : \omom \to 2^\omega$ via
\[ \varphi (f) (n) : = f(n)\; \mod\; 2 \]
for $f\in\omom$ and $n\in\omega$. Note that $\varphi$ is onto, 
continuous and open, regardless of whether we topologize 
$\omom$ with ${\cal D}$ or ${\cal B}$. (Of course, $2^\omega$
always carries the topology ${\cal C}$).

Now let $A \subseteq 2^\omega$ be a ${\cal C}$--nonmeager 
set in $\Gamma$. It suffices to show that
there is $s \in 2^{<\omega}$ such that $[s] \cap A$ is comeager
in $[s]$. Since $\varphi$ is continuous when going from ${\cal B}$
to ${\cal C}$ and $\Gamma$ is topologically reasonable,
$B = \varphi^{-1} (A)$ is in $\Gamma$ as well. As $\varphi$ is onto,
continuous and open as a map from $\la \omom, {\cal D}\ra$ to
$\la 2^\omega , {\cal C}\ra$, $B$ is ${\cal D}$--nonmeager. 
By assumption we find $[N,f]$ such that $[N,f] \cap B$ is comeager in 
$[N,f]$. Hence there is a $G_\delta$--set $C \subseteq [N,f] \cap B$
dense in $[N,f]$. Assume $C = \bigcap_i C_i$ where the $C_i$
form a decreasing sequence of open sets, and $C_i = \bigcup_{
|\sigma| = i} [N_\sigma , f_\sigma]$ (where $N_{\la\ra} = N$
and $f_{\la\ra} = f$) are such that 
\begin{enumerate}
\item \label{DiiDdense} $\bigcup_j [N_{\sigma \et \la j\ra} , f_{\sigma 
\et\la j\ra}] 
\subseteq [N_\sigma , f_\sigma]$ is dense for all $\sigma$, and
\item \label{DiiDgeq} $N_\sigma \geq i$ for all $i$ and $\sigma$ 
with $|\sigma | = i$.
\end{enumerate}
It is clear that all the $C_i$ can be written in this form.

Let $\bar\varphi : \omlom \to 2^{<\omega}$ be defined by
\[ \bar\varphi (s) (n) : = s (n) \;\mod\; 2 \]
for $s\in\omlom$ and $n < |s|$, and put $s_\sigma = \bar\varphi
(f_\sigma | N_\sigma )$ and $s = s_{\la\ra}$. Next find
$H_\sigma \subseteq \omega$ such that 
\begin{enumerate} \setcounter{enumi}{2}
\item \label{DiiCpd} the $[s_{\sigma\et \la j\ra} ]$ for 
$j \in H_\sigma$
are pairwise disjoint, and
\item \label{DiiCdense} $\bigcup_{j\in H_\sigma} [s_{\sigma \et\la
j\ra}]$
is dense in $[s_\sigma]$.
\end{enumerate}
Again this is easily done by \ref{DiiDdense} above.
Now define recursively $E_0 = \{ \la\ra \}$, $E_{i+1} = \bigcup_{\sigma
\in E_i} \{ \sigma \et \la j\ra : \;  j \in H_\sigma\}$, put 
$D_i = \bigcup_{\sigma \in E_i} [s_\sigma]$, and let 
$D = \bigcap_i D_i$. By \ref{DiiCdense}, $D$ is dense in
$[s]$. We claim that $D \subseteq A$, thus completing the proof.

Given $x\in D$, there is a unique $y\in\omom$ such that $x \in
[s_{y | n}]$ for all $n\in\omega$, by clause \ref{DiiCpd}.
Furthermore, \ref{DiiDgeq} entails that $\bigcap_n [s_{y|n}]
= \{ x \}$. Now, $\bigcap_n [N_{y|n}, f_{y|n} ]$ also
contains a unique element $g\in\omom$, by \ref{DiiDdense}
and \ref{DiiDgeq}. Clearly, $\varphi (g) = x$.
Since $g \in C \subseteq B = \varphi^{-1} (A)$, we get
$x \in A$, as required. \qed
Note that this result is nothing but a topological version of
the well--known fact that if $f\in\omom$ is {\sc Hechler}
over a model $\frnM$ of set theory, then $\varphi (f)$
is {\sc Cohen} over $\frnM$. For the next result (Theorem \ref{DiiDDiiL}
below), we need the following notion from \cite{BHS95} (p. 294):
\begin{sdef}
Let $\bar W = \la w_\sigma , s_\sigma : \; \sigma \in \omlom \ra$
be such that
\begin{itemize}
\item $\dom(s_{\la \ra})$ and $w_\sigma$  are finite subsets of $\omega$
\item $s_\sigma : w_{\sigma | |\sigma| - 1} \to\omega$ for
$\sigma \neq \la\ra$, $s_{\la\ra} : \dom (s_{\la\ra} )\to \omega$ are
functions
\item $\omega = \dom (s_{\la\ra}) \cup \bigcup_{n\in\omega} w_{f|n}
$ for all $f \in \omom$, the union being pairwise disjoint
\item $s_\sigma (i) > \sigma ( |\sigma| - 1)$ for all $i \in w_{\sigma
||\sigma | -1}$ and all $\sigma$
\end{itemize}
Then we can define the set $C = C(\bar W) \subseteq \omom$ such that
$g \in C $ iff $g= \bigcup_n s_{f|n}$  for some $f\in\omom$.
$C$ is called a {\it nice set}; it is necessarily closed and dominating.
\end{sdef}
\begin{sthm}[\cite{BHS95}{\rm , Theorem 1.1}] \label{nicethm}
Every dominating analytic set contains a nice set.
\end{sthm}
\begin{sthm} \label{DiiDDiiL}
For any topologically reasonable pointclass $\Gamma$,
$\Gamma(\DDD)$ implies $\Gamma(\LLL)$.
\end{sthm}
{\bf Proof :} \\
Let $A\in \Gamma$ and let $T$ be
a {\sc Laver} tree. By \ref{weakM}, we can assume  $T = \omega^{
< \omega}$.
We have to find a {\sc Laver} tree $S \leq T$ such
that either $[S] \subseteq A$ or $[S] \cap A = \emptyset$.
We define a function $\varphi :\omom \to \omom$
recursively by
\[ \begin{array}{ccc}
\varphi (x) (0) & = & x (0) \\
\varphi (x) (n+1) & = & x(\varphi (x) (n))
\end{array} \]
Clearly $\varphi$ is continuous. Put $B : = \varphi^{-1} (A)$.
By assumption $B \in \Gamma$. Hence $B$ has the property of 
{\sc Baire} in the topology ${\cal D}$. Thus we can find an open set
$[N,f]$ in ${\cal D}$ such that either $B$ or $\omom \setminus
B$ is ${\cal D}$--comeager in $[N,f]$. Without loss the former
holds. Hence there is a $G_\delta$--set $C \subseteq [N,f]
\cap B$ dense in $[N,f]$. Note that $C$ must be dominating in $\omom$,
for otherwise we could find $g\in\omom$ above $f$ with
$[N,g] \cap C = \emptyset$, contradicting $C$'s density. By
Theorem \ref{nicethm}, $C$ contains a nice set $D = C(\bar W)$
where $\bar W = \la w_\sigma , s_\sigma : \; \sigma \in \omlom \ra$.
Since $D \subseteq B$, we get $\varphi [D] \subseteq A$.
We are left with showing that $\varphi[D]$ contains the set of branches through
a {\sc Laver} tree $S$.

To this end, define recursively a function $\bar\varphi$ with
range $\omlom$ and domain all functions from a finite subset of $\omega$ to
$\omega$, as follows. If $0 \notin \dom (s)$, let $\bar\varphi (s) = \la\ra$.
Otherwise put
$$\bar\varphi (s) (0) : = s(0)$$
Assume $\bar\varphi (s) (i)$ has been defined. If $\bar\varphi (s)
(i) \notin \dom (s)$, we're done and have $|\bar\varphi (s) | = i+1$.
Otherwise, put
$$\bar\varphi (s) (i+1) : = s (\bar\varphi (s) (i))$$
Now construct recursively a {\sc Laver} tree $S$ such that for any
$t \in S$ there is $\sigma\in\omlom$ such that
$\bar\varphi (\bigcup_{j \leq |\sigma|} s_{\sigma | j} ) = t$
$(\star)$. Clearly $(\star)$ implies $[S] \subseteq \varphi [D]$.

First put $t : = \bar\varphi (s_{\la\ra})$ into $S$.
Then assume $t \in S$ has property $(\star)$ with witness $\sigma$.
We have to define the successors of $t$ in $S$.
Put $s := \bigcup_{j \leq |\sigma|} s_{\sigma | j}$.
Then by definition of $\bar\varphi$, $t(|t| - 1) \notin
\dom (s)$. Hence there is $\tau \supseteq \sigma$ minimal such that $t(|t|
- 1) \in w_\tau$. Now, if $n$ is large enough, we will
have $\bar\varphi (s_n) (|t|) = s_n (t (|t| - 1)) \notin \dom (s_n)$
where $s_n = \bigcup_{j \leq |\tau| + 1} s_{\tau\et\la n \ra | j}$.
Therefore $t_n = \bar\varphi (s_n)$ for such $n$ 
has length $|t| + 1$ and also has property $(\star)$ with witness $\tau\et\la n
\ra$. Thus we can put such $t_n$ into $S$. This completes the recursive
construction of the {\sc Laver} tree $S$, and the proof of the Theorem.
\qed
Results like \ref{DiiDDiiC} and \ref{DiiDDiiL} can be subsumed
in the following diagram.
\begin{scor} \label{Gamma}
Let $\Gamma$ be a topologically reasonable pointclass. Then one
has the following implications:
\[ \begin{array}{ccccc}
\Gamma(\DDD) & \Longrightarrow & \Gamma(\CCC) && \\
\Downarrow && \Downarrow &&\\
\Gamma(\LLL) & \Longrightarrow & \Gamma (\MMM) & \Longrightarrow & 
\Gamma(\SSS) 
\end{array} \]
\end{scor}
{\bf Proof :} \\
$\Gamma (\DDD) \Longrightarrow \Gamma (\CCC)$ and
$\Gamma (\DDD) \Longrightarrow \Gamma (\LLL)$ were proved
in Theorems \ref{DiiDDiiC} and \ref{DiiDDiiL}, respectively.

The directions
$\Gamma(\LLL)  \Longrightarrow  \Gamma (\MMM) 
\Longrightarrow  \Gamma(\SSS)$ are easy consequences of \ref{weakM}.
To see e.g. the second implication, let $A \subseteq 2^\omega$
be a set in $\Gamma$, and let $S \subseteq 2^{< \omega}$ be a
{\sc Sacks} tree. By \ref{weakM}, we can assume $S = 2^{< \omega}$. Let
$\varphi : \omom \to 2^\omega$ be the canonical map which identifies
the {\sc Baire} space with the irrationals in $2^\omega$
({\it i.e.}, the $x \in 2^\omega$ such that $\{ i : \; x(i) = 1 \}$ is
infinite). Since $\varphi$ is continuous, $\varphi^{-1} (A)$ belongs
to $\Gamma$. Hence we can find $M \in \MMM$ with $[M] \subseteq
\varphi^{-1} (A)$ or $[M] \subseteq \omom \setminus \varphi^{-1}
(A)$. Assume without loss the former. Since $\varphi [M]$ is an uncountable 
$G_\delta$--set, we can find $T \in \SSS$ with $[T] \subseteq \varphi [M]
\subseteq A$, as required.

 To see that $\Gamma (
\CCC)$ implies $\Gamma (\MMM)$, simply note that every non--meager
set with the property of {\sc Baire} contains the set of branches through 
a superperfect tree.\qed
We sketch another connection between two regularity properties which
we shall need in section \ref{Laversec} when dealing with {\sc Laver} forcing.
To this end,  we introduce the following 
three notions the first of which is Definition 2.1 in \cite{GRSS} while 
the last is on p. 296 in \cite{BHS95}: 
\begin{sdef}
\begin{enumerate} 
\item A set $A \subseteq \omom$ is called {\it strongly dominating} iff
\[ \forall f \in\omom \; \exists x \in A \; \forall^\infty
k : f(x(k-1) ) < x(k) \]
\item A set $A \subseteq \omom$ is called {\it $\ell$--regular} if either
$A$ contains the set of branches through a Laver tree or $A$ is not
strongly dominating.
\item A set $A \subseteq \omom$ is called {\it strongly $u$--regular}
if either $A$ contains a nice set or $A$ is not dominating. 
\end{enumerate}
\end{sdef}
Note that the second and third notions are very similar, and analogous facts
can be proved about both. It was shown in Lemma 2.3 of \cite{GRSS} that every
{\sc Borel} set is $\ell$--regular. Standard modifications of the
game--theoretic argument used in the proof ({\sc Solovay}'s unfolding trick) 
show the same conclusion is true for analytic sets --- this is, of course,
analogous to Theorem \ref{nicethm} above, but it's also a consequence
of \ref{nicethm} and the following proposition:
\begin{sprop} \label{regprop}
For a topologically reasonable pointclass $\Gamma$, strong $u$--regularity
for $\Gamma$ implies $\ell$--regularity for $\Gamma$.
\end{sprop}
{\bf Proof :} \\
Let $A \in \Gamma$ be strongly dominating. Let $\varphi :
\omom\to\omom$ be the function constructed in the proof of Theorem
\ref{DiiDDiiL}. By the proof of \ref{DiiDDiiL}, it suffices to show
that $B : = \varphi^{-1} (A)$ is dominating --- for then we can use
strong $u$--regularity to get a nice set $C \subseteq B$ and
the argument of \ref{DiiDDiiL} shows that
$\varphi [C] \subseteq A$ contains a {\sc Laver} tree.

To see that $B$ is dominating, let $g \in\omom$ be an arbitrary increasing
function such that $\varphi (g) (j - 1) \geq j$ for all $j$.
Find $x \in A$ such that $x (n+1) > \varphi (g) (x(n))$ for all
$n\in\omega$. Define $y \in\omom$ such that
\[ \begin{array}{ccccc}
y (0) & = & x (0) &&\\
y(i) & =  & x(1) & \hspace{ 2cm} & \mbox{ for } 1 \leq i \leq x(0) \\
y(i) & = & x(n+1) &\hspace{ 2cm} & \mbox{ for } n \geq 1 \mbox{ and }
x(n-1) < i \leq x(n) \\
\end{array} \]
Then $\varphi (y) = x$ and hence $y \in B$. Furthermore,
\[y(i) = x(n+1) > \varphi (g) (x(n)) = g( \varphi (g) (x(n) - 1)) \geq
g(x(n)) \geq g(i) \]
for $x(n-1) < i \leq x(n)$, because $\varphi (g) (x(n) - 1)
\geq x(n)$. Thus we have $y \geq^* g$, as required. \qed

%%%%%%%%%%%%%%%%%%%%%%%%%   Section 4   %%%%%%%%%%%%%%%%%%%%%%%%%%

\section{{\sc Laver} Measurability} \label{Laversec}
In contrast to the topological forcings (see section \ref{Hechsec}), 
for the three non--topological forcings the notions of
$\Db{1}{2}$-- and $\Sb{1}{2}$--measurability are equivalent. For {\sc Laver}
forcing we will prove:
\begin{sthm} \label{Laverchar}
The following are equivalent:
\begin{enumerate}
\item $\forall a\in\omom : \omom\cap{\bf L}[a]$ is $\sigma$--bounded
in $\omom$
\item $\DiiL$
\item $\SiiL$
\end{enumerate}
\end{sthm}
For our proof, we need the following characterization part of which
is a consequence of \ref{regprop}.
\begin{sprop} \label{ellreg}
The following are equivalent:
\begin{enumerate}
\item Every $\Sb{1}{2}$--set is strongly $u$--regular
\item Every $\Sb{1}{2}$--set is $\ell$--regular
\item \label{sigmabound} $\forall a \in\omom: \omom \cap {\bf L} [a]$ is
$\sigma$--bounded in $\omom$
\end{enumerate}
\end{sprop}
{\bf Proof :}\\
(i)$\Rightarrow$(ii): By Proposition \ref{regprop}.\\
(ii)$\Rightarrow$(iii): This will follow from (ii)$\Rightarrow$(i)
in \ref{Laverchar}, because $\Sb{1}{2} -\ell$--regularity clearly
implies $w\SiiL$, and hence $\SiiL$ by Lemma \ref{weakM}.\\
(iii)$\Rightarrow$(i): This was remarked on p. 296 in
\cite{BHS95}. The proof is identical to the proof of Theorem 4.2 of \cite{S94}.
\qed
{\bf Proof of \ref{Laverchar} :}\\
(i)$\Rightarrow$(iii): This is immediate from the direction
(iii)$\Rightarrow$(ii) in \ref{ellreg}.\\
(ii)$\Rightarrow$(i): 
Suppose we had an $a$ so that ${\bf L}[a]\cap\omom$ is not $\sigma$--bounded, 
then:
\[\forall x\in\omom\;\exists y\in{\bf L}[a]\cap\omom\;\exists^\infty
n\in\omega:
y(n)>x(n)\]
Let $\la g_\alpha :\alpha<\omega_1\ra$ be the 
$\Sl{1}{2}(a)$--good well--ordering of ${\bf L}[a]$. 
From this we can define a $\Sl{1}{2}(a)$--scale 
$\la f_\alpha :  \alpha < \omega_1 \ra$ in 
${\bf L}[a]$\footnote{{\it I.e.} a dominating
subset of $\omom\cap{\bf L}[a]$ well--ordered by $\leq^*$.} which is 
unbounded in $\omom$ and additionally has
the property
$$\forall\alpha<\omega_1\;\forall n<\omega\;:\;
f_{\alpha+1} (n) \geq  f_\alpha (n+1),$$
by standard tricks.
With this scale of reals we define the following sets:
\begin{sdef}
$$x\in A_\alpha:\iff (\forall\beta<\alpha:x^*\geq f_\beta)\wedge
\exists^\infty n (x(n)<f_\alpha(n))$$
$$A:=\bigcup_{\alpha\mbox{{\footnotesize{} is even}}}A_\alpha$$
$$B:=\bigcup_{\alpha\mbox{{\footnotesize{} is odd}}}A_\alpha$$
As usual, limit ordinals are counted as even.
\end{sdef}
As is easily checked, the family $\la A_\alpha : \alpha<\omega_1\ra$
is pairwise disjoint and covers all of $\omom$. Therefore $A$ and $B$ are
complementary. 
Because the scale was $\Sl{1}{2}(a)$, both $A$ and $B$ 
are $\Dl{1}{2}(a)$--sets.

Next take a {\sc Laver} tree $L$. Without loss of generality we may assume
that for all nodes $s\in L$ we have the following property:
$$\forall t\in\Succ(s) : t(|s|) > s(|s|-1)$$
Now we define recursively for any $s\in L$:
$$g_s(n) := s(n)\mbox{ for }n<|s|$$
$$g_s(n) := \min\{t(n) : t\in\Succ(g_s|n)\}\mbox{ for }n\geq|s|$$
Then $g_s\in[L]$ and $g_s|m\in L$ for all $m<\omega$. Because of our
assumption on $L$ the $g_s$ are strictly increasing after the stem of
$L$.

Now find $\alpha$ so that $ f_{\alpha}$ lies infinitely
often above each $g_s$ for $s\in L$. This is possible by the
unboundedness of the sequence $\la f_\alpha : \alpha < \omega_1 \ra$.
To prove the theorem
we have to show that for the arbitrarily chosen {\sc Laver} tree $L$
there is a branch through $L$ in $A$ as well as in $B$. To this end
we will prove the following stronger claim:
\begin{quote}
Let $\gamma\geq\alpha$. Then there is $x\in [L]\cap A_{\gamma + 1}$.
\end{quote}

For this, we make the following recursive construction.
Define $s_0$ to be the stem of $L$. If $s_i$ is already defined, 
choose $t\in\Succ(s_i)$
so that $t(|s_i|)\geq  f_\gamma(|s_i|)$. 
We know that $ f_{\gamma+1}$ has infinitely many points where it
is above $g_t$. Take $n\geq|s_i|$ minimal with this property. Then 
for all $m$ with
$|s_i|<m\leq n$:
$$f_\gamma (m)\leq  f_{\gamma+1} (m-1) \leq g_t(m-1)< g_t(m)
$$ Since we also have
$$ f_\gamma (|s_i|) \leq t(|s_i|) = g_t(|s_i|)$$
we know that $f_\gamma$ lies below $g_t$ between $|s_i|$ and $n$ and
$ f_{\gamma+1}(n)>g_t(n)$. Hence define $s_{i+1}:=
g_t|n+1$.

Now we put $x:=\bigcup_{i\in\omega} s_i$. According to the
construction, $x$ 
dominates $f_\gamma$ and $x(|s_i|-1)<
f_{\gamma+1}(|s_i|-1)$ for all $i<\omega$. Thus $x\in A_{\gamma +1}$.
Because all $s_i$ were in
$L$, we have $x\in[L]$.

Since $\gamma\geq\alpha$ was arbitrary we have elements of $[L]$
both in $A$ and $B$, whence $A$ and $B$ cannot be $\LLL$--measurable.\qed

%%%%%%%%%%%%%%%%%%%%%   section 5    %%%%%%%%%%%%%%%%%%%%%%%%%%%%%%%

\section{{\sc Hechler}-Forcing} \label{Hechsec}

This section is devoted to proving the characterizations
of $\SiiD$ and $\DiiD$ mentioned in the Introduction.
For this we will need the well--known characterizations
for $\SiiC$ and $\DiiC$:
\begin{sthm}[{\sc Solovay}]\label{SoloSb12}
The following are equivalent:
\begin{enumerate}
\item $\SiiC$
\item $\forall a\in\omom:\NC{a}$ is meager
\item $\forall a\in\omom:\CohL{a}$ is comeager
\end{enumerate}
\end{sthm}
\begin{sthm}[\cite{JSh89}{\rm , Theorem 3.1}]\label{DiiC}
The following are equivalent:
\begin{enumerate}
\item $\DiiC$
\item $\forall a\in\omom: \CohL{a}\not=\emptyset$
\end{enumerate}
\end{sthm}
For  proofs {\it cf.} \cite{BJ95}, p. 457 and p. 452{\it sqq.}
To get from these equivalences  results about $\DiiD$ and $\SiiD$ we need a
connection between $\CCC$ and $\DDD$. This connection is provided
by the theorems of {\sc Miller} and {\sc Truss}:
\begin{sthm}[\cite{T77}{\rm , Lemma 6.2}]\label{Truss}
If $\frM$ is a ZFC--model, $d$ a dominating real over
$\omom\cap\frM$ and $c\in\omom$ {\sc Cohen} over
$\frM[d]$, then $c+d$ is a {\sc Hechler} real over $\frM$.
\end{sthm}
\begin{sthm}[\cite{T77}{\rm , Theorem 6.5}]\label{Truss2}
Let $c$ be a {\sc Cohen} real over $\frM$ and $d$ dominating
over $\frM[c]$. Then the set of all {\sc Cohen} reals
over $\frM$ is comeager.
\end{sthm}
\begin{sthm}[\cite{M81}{\rm  , Theorem 1.2}]\label{Miller}
Consider the partial orderings $\la\omom,\leq^*\ra$ and
$\la (c^0),\subseteq\ra$. 
For $f\in\omom$ define
$$\tilde f(i):=\max\{f(j)+1:j\leq i\}$$
Then the function $T:\omom\to (c^0) :
f\mapsto\{x\in\omom : x\leq^*\tilde f\}$ satisfies:\\
For every set $X$ bounded in $(c^0)$,
$T^{-1}(X)$ is bounded in $\omom$.
\end{sthm}
For a proof {\it cf.} \cite{BJ95}, p. 39{\it sq}.
As an easy corollary we get:
\begin{scor}\label{SiiCDom}
Suppose that for every $a\in\omom$ there is a {\sc Cohen} real over
${\bf L}[a]$.
Then the following are equivalent:
\begin{enumerate}
\item $\SiiC$
\item $\forall a\in\omom$ : ${\bf L}[a]\cap\omom$ is $\sigma$--bounded
\end{enumerate}
\end{scor}
{\bf Proof :}\\
(i)$\Rightarrow$(ii): According to \ref{SoloSb12} $\NC{a}$ is meager,
therefore ${\cal M}:=\la A_c:c\in\BCC{a}\ra$ is a bounded family in $(c^0)$.
Hence according to \ref{Miller} $T^{-1}({\cal M})$ a bounded family in
$\omom$. Because we constructed $T$ in ${\bf L}[a]$, the set $T(x)$ is
a meager set coded in ${\bf L}[a]$ for any $x\in\omom\cap {\bf L}[a]$.
So we have $\omom\cap {\bf L}[a]\subseteq T^{-1}({\cal M})$, and
therefore the real numbers of ${\bf L}[a]$ are bounded.\\[0.1cm]
(ii)$\Rightarrow$(i): Following the assumption we have over each
${\bf L}[a]$ a {\sc Cohen} real $c_a$ and a dominating real
over ${\bf L}[a][c_a]$, hence we have with \ref{Truss2}:
The set of all {\sc Cohen} reals over ${\bf L}[a]$ is comeager. 
Because $a$ was arbitrary, the claim follows from \ref{SoloSb12}.\qed
The following result which is a consequence of earlier theorems is the
cornerstone of the proof of Theorem \ref{DiiD} below.
\begin{scor}\label{DiiDbeschr}
$\DiiD\Rightarrow\forall a\in\omom:\omom\cap 
{\bf L}[a]$ is $\sigma$--bounded
\end{scor}
{\bf Proof :}\\
Follows from Theorems \ref{DiiDDiiL} and \ref{Laverchar}. \qed
Of course, this result can also be proved directly, without
any reference to {\sc Laver} forcing.

\begin{sthm}\label{DiiD}
The following are equivalent:
\begin{enumerate}
\item $\DiiD$
\item $\forall a\in\omom:\HechL{a}\neq\emptyset$
\item $\SiiC$
\end{enumerate}
\end{sthm}
{\bf Proof :}\\
(i)$\Rightarrow$(iii):
Because of \ref{DiiDDiiC} we have $\DiiC$,
especially there is a {\sc Cohen} real over each ${\bf L}[a]$
according to \ref{DiiC}. 
Because of \ref{DiiDbeschr} we get a dominating real over
${\bf L}[a]$, and from the {\sc Cohen} and the dominating real
we conclude $\SiiC$
via \ref{SiiCDom}.\\
(iii)$\Rightarrow$(ii): According to \ref{SiiCDom}, $\SiiC$ implies 
the existence of a dominating real over each ${\bf L}[a]$.
Together with the {\sc Cohen} real we get from
\ref{SoloSb12}, we get with \ref{Truss} a 
{\sc Hechler} real over each ${\bf L}[a]$.\\
(ii)$\Rightarrow$(i):
This is exactly the same proof as in the corresponding direction of
\ref{DiiC} ({\it cf.} \cite{JSh89}, Theorem 3.1,
or \cite{BJ95}, p. 452{\it sqq}.).
\qed
Notice that this result can be looked at as a ``projective" version
of the combinatorial result that the covering number of the
ideal $(d^0)$ is equal to the additivity of $(c^0)$ ({\it cf.}
\cite{LR95}, Theorem 3.6). We are now heading towards
a characterization of $\SiiD$. Apart from what has been proved so
far, the following combinatorial tool is essential. 
Let ${\cal A}$ be an almost disjoint system of subsets
of $\omega$. We define for
$A\in{\cal A}$:
$$X_A :=\{x\in\omom : \ran(x)\cap A =\emptyset\}$$
As one can easily see $X_A$ is a closed nowhere dense
set in ${\cal D}$.
\begin{sthm}[\cite{LR95}{\rm , Theorem 6.2}]\label{LRThm}
If $X$ is a $\DDD$--null set, then there are at most countably many 
$A\in{\cal A}$, so that $X_A\subseteq X$.
\end{sthm}
\begin{slem}\label{Trichotomie}
Suppose that the set $\HechL{a}$ has the {\sc Baire} property,
then it is either meager or comeager in ${\cal D}$.
\end{slem}
{\bf Proof :}\\
Suppose $\HechL{a}$ has the {\sc Baire} property in ${\cal D}$.
If $\HechL{a}$ is not meager then there is an open set $[N,f]$ in which
$\HechL{a}$ is comeager. It suffices to show that below  each open set
$[M,g]$ there is another open set in which 
$\HechL{a}$ is comeager. Take an open set $[M,g]$ and 
define $\tilde f(i):=g(i)$ for $i<M$ and $\tilde f(i+M):=
f(i+N)$.
Then --- since changing finite initial segments does not change the
property of being {\sc Hechler} --- $\HechL{a}$ is still comeager
in $[M,\tilde f]$ and therefore in every subset.
Obviously one can find in $[M,\tilde f]$ a real $h$ lying completely
above $g$. Then $[M,h]$ is a subset of both $[M,g]$ and 
$[M,\tilde f]$.
Hence $\HechL{a}$ is comeager in $[M,h]$, as required.\qed

\begin{sthm} \label{HechSigma}
The following are equivalent:\begin{enumerate}
\item $\SiiD$
\item $\forall a\in\omom : \aleph_1^{{\bf L}[a]} <\aleph_1$
\end{enumerate}
\end{sthm}
We divide the proof into two parts, one using 
\ref{DiiD} and the other using \ref{LRThm}.

\begin{sprop}
\[\SiiD\iff\forall a\in\omom : \HechL{a}\in (d^1)\]
\end{sprop}
{\bf Proof :}\\
``$\Leftarrow$": Let $A$ be a $\Sl{1}{2}(a)$--set. By
a well--known result of \cite{S70} ({\it cf.} also \cite{Jech},
p. 545), there is a {\sc Borel} set $B$ such that
$A \cap \HechL{a} = B \cap \HechL{a}$. Since $\HechL{a}$ is
${\cal D}$--comeager, it follows that $A$ has the {\sc Baire} property in
${\cal D}$.\\
``$\Rightarrow$": The set $\ND{a}$
is a $\Sl{1}{2}(a)$--set, because
$$x\in\ND{a}\iff \exists c\in\BCD{a}:
x\in A_c$$
Hence it is $\DDD$--measurable. According to \ref{Trichotomie}
it is either comeager or meager. We have to exclude the case that
it is comeager. Suppose $\HechL{a} = \omom \backslash
\ND{a}$ is a meager set. Then it is 
included in some meager set coded in some
${\bf L}[b]$, hence there are no {\sc Hechler} reals over ${\bf L}[a,b]$
anymore.
But $\SiiD$ implies by \ref{DiiD} the existence of a {\sc Hechler} real over
${\bf L}[a,b]$, a contradiction.\qed

\begin{sprop}
\[\forall a\in\omom : \HechL{a}\in (d^1)\iff
\forall a\in\omom : \aleph_1^{{\bf L}[a]} <\aleph_1\]
\end{sprop}
{\bf Proof :}\\
``$\Leftarrow$": If
$(2^{\aleph_0})^{{\bf L}[a]}=(\aleph_1)^{{\bf L}[a]}$ is countable,
then there are at most countably many codes for $\DDD$--null sets
in ${\bf L}[a]$, hence $\ND{a}$ is a $\DDD$--null set.\\
``$\Rightarrow$": Suppose $(\aleph_1)^{{\bf L}[a]}=\aleph_1$ 
for some $a$.
Then there is in ${\bf L}[a]$ an almost disjoint family $\cal A$ with
$|{\cal A}| = (2^\omega)^{{\bf L}[a]} = (\aleph_1)^{{\bf L}[a]}=\aleph_1$.
We know that the sets $X_A$ are {\sc Hechler}-null sets
in ${\bf L}[a]$.
Because of that $\ND{a}$ contains all of the $X_A$ and hence more
than countably many of these sets.
\ref{LRThm} shows that $\ND{a}$ is not {\sc Hechler}--null.\qed

%%%%%%%%%%%%%%%%%%%%%%%%%%   Section 6   %%%%%%%%%%%%%%%%%%%%%%%%%%%%%

\section{{\sc Miller} Measurability} \label{Millersec}
The main goal of this section is the proof of the following
characterization:
\begin{sthm}\label{MCharac}
The following are equivalent:
\begin{enumerate}
\item $\forall a\in\omom :$ $\omom\cap{\bf L}[a]$ is not dominating in 
$\omom$
\item $\DiiM$
\item $\SiiM$
\end{enumerate}
\end{sthm}
From \cite{BHS95} we introduce the following 
notion:
\begin{sdef}
A set $B\subseteq\omom$ is called {\it $w$--regular} 
if either $B$ contains the set
of branches through a superperfect tree or $B$ is not dominating.
\end{sdef}
An old result from \cite{K77} ({\it cf.} Theorem 4) yields
({\it cf.} Proposition 2.3 in \cite{BHS95}):
\begin{sthm}[{\sc Kechris}/{\sc Spinas}] \label{wreg}
The following are equivalent:
\begin{enumerate}
\item Every $\Sb{1}{2}$ set is $w$--regular
\item $\forall a\in\omom :$ $\omom\cap{\bf L}[a]$ is not dominating in 
$\omom$
\end{enumerate}
\end{sthm}
{\bf Proof of \ref{MCharac} :}\\
(i)$\Rightarrow$(iii): With Theorem \ref{wreg} we immediately get 
$w\SiiM$ and with that via \ref{weakM} $\SiiM$.\\
(ii)$\Rightarrow$(i):
Let $\la\sigma_n : n<\omega\ra$ be an enumeration of $\omlom$.
Let $\code : \omlom\to\omega$ be defined by
$$\code(\sigma) = n \iff \sigma = \sigma_n$$
With a given superperfect tree $T\subseteq \omlom$ we associate 
a function $f_T\in\omega^{(\omlom)}$ and a sequence
$\la\tau^T_\sigma : \sigma\in\omlom\ra$ of elements of $T$ using
the following recursion:
$$f_T(\la\ra) := \code(\stem(T))$$
$$\tau^T_{\la\ra} := \stem (T)$$
\begin{eqnarray*}
f_T(\sigma) & := & \min\{ n : \stem(T)\et\tau^T_{\la\sigma(0)\ra}\et\dots
\et\tau^T_{\sigma|(|\sigma|-1)}\et\sigma_n\in\Split(T)\\
&& \mbox{ and }
\sigma_n(0)>\sigma(|\sigma|-1)\}
\end{eqnarray*}
$$\tau^T_\sigma := \sigma_{f_T(\sigma)}$$
Call $f\in \omega^{(\omega^{<\omega})}$ {\it fast} iff for all $\sigma \in
\omlom$ there is $\tau\in\omlom$ with $\tau (0) > \sigma (
|\sigma| - 1)$ and $\code (\tau) < f (\sigma)$.
Given $f\in\omega^{(\omlom)}$ fast, $g\in\omom$ and a natural number
$m\in\omega$, we  define recursively the tree $T=T(f,g,m)$:
$$T_0 := \{\sigma_i|\ell : i<m, \ell\leq|\sigma_i|\}$$
$$\tilde g(0) := m$$
$$T_1 := \{p\et(\sigma_i|\ell) : p\in T_0, i< f(\la\tilde g(0)\ra),
\ell\leq|\sigma_i|, \sigma_i(0)>\tilde g(0)\}$$
Note that $T_1\setminus T_0\neq\emptyset$. In the $n$th
step we put:
$$\tilde g(n)  := \max\{g(j) : j\leq{\rm height}(T_n)\}$$
$$T_{n+1} := \{p\et(\sigma_i|\ell) : p\in T_n\setminus T_{n-1}, 
i< f(\la\tilde g(0),\dots, \tilde g(n)\ra),
\ell\leq|\sigma_i|, \sigma_i(0)>\tilde g(n)\}$$
Let $T=\bigcup_{n\in\omega} T_n$. Notice
that $T$ is a finitely branching tree, that is $[T]$ is compact.
Obviously, no branch of $T(f,g,m)$ is eventually dominated by $g$.

\begin{slem}\label{MLem}
If $S$ is a superperfect tree, $f$ is fast and $f_S<^* f$, then there is
an $m\in\omega$ such that $[S]\cap[T(f,g,m)]\neq\emptyset$.
\end{slem}
{\bf Proof :}\\
Choose $m\in\omega$ such that $f_S(\la\ra) < m$ and 
$f_S(\la m\ra\et\sigma)< f(\la m\ra\et\sigma)$ for all $\sigma$. We construct
a branch belonging to both trees as follows:
$$\tau_0 := \stem(S)$$
$$\tau_1 := \tau^S_{\la m\ra}$$
$$\tau_n := \tau^S_{\la m,\tilde g(1),\dots,\tilde g(n-1)\ra},$$
where $\tilde g$ is constructed from $g$ as above. Then $\tau_0\in T_0$
(by $f_S(\la\ra)<m$), $\tau_0\et\tau_1\in T_1$ (by $f(\la m\ra)>f_S(\la m\ra)$
and $\tau_1(0)=\tau^S_{\la m\ra}(0)>m$), and so on.\\
Thus $x:=\tau_0\et\tau_1\et\dots\et\tau_n\et\dots\in[S]\cap[T(f,g,m)]$.\qed
We are now ready to complete the proof of \ref{MCharac}.

Suppose the reals of ${\bf L}[a]$ were dominating. Let $\la g_\alpha:
\alpha<\omega_1\ra$ and $\la g_\alpha':
\alpha<\omega_1\ra$ be the $\Sl{1}{2}(a)$--enumerations of
$\omom\cap{\bf L}[a]$ and $\omega^{(\omlom)}\cap{\bf L}[a]$, respectively.
We construct recursively $\la f_\alpha : \alpha<\omega_1\ra\subseteq
\omom\cap{\bf L}[a]$ and auxiliary 
$\la h_\alpha : \alpha<\omega_1\ra\subseteq
\omega^{(\omlom)}\cap{\bf L}[a]$, such that for $\alpha<\beta$:
\begin{enumerate}
\item $h_\beta$ eventually dominates $g_\beta'$ and $h_\alpha$
\item $f_\alpha$ is fast
\item $f_{\alpha+1}$ eventually dominates all branches of all trees
$T(h_\alpha, f_\alpha, m)$ for $m<\omega$ (possible by compactness)
\item $f_\beta$ eventually dominates $f_\alpha$ and $g_\beta$
\item $f_\alpha$ and $h_\alpha$ are $<_{{\bf L}[a]}$--minimal with these
properties
\end{enumerate}
We form
$$A:=\{y\in\omom : \min\{\alpha : y<^* f_\alpha\}\mbox{ is even}\}$$
$$B:=\{y\in\omom : \min\{\alpha : y<^* f_\alpha\}\mbox{ is odd}\}$$
Then both $A$ and $B$ have $\Sl{1}{2}(a)$--definitions, $A\cap B=\emptyset$
and $A\cup B=\omom$, because we worked through the $\la g_\alpha : \alpha<
\omega_1\ra$.

Since the reals in ${\bf L}[a]$ are dominating, $A\cup B=\omom$
still holds in the real world ${\bf V}$. Thus $A$ and $B$ are both
$\Dl{1}{2}(a)$ in ${\bf V}$. We now show that $A$ and $B$ are
super--{\sc Bernstein}--sets, {\it i.e.} for each superperfect
tree $T$ we have $A\cap [T]\neq\emptyset\neq B\cap [T]$. For this
purpose let
$T\in{\bf V}$ be superperfect.
There is an $\alpha<\omega_1$ such that $f_T<^* h_\beta$ for all
$\beta\geq\alpha$, because $\omom\cap{\bf L}[a]$ is dominating and we worked
through the $\la g_\alpha' : \alpha\in\omega_1\ra$.

Thus there are (by \ref{MLem}) $m$ and $m'\in\omega$ with
$$[T]\cap[T(h_\alpha, f_\alpha, m)]\neq\emptyset$$
and
$$[T]\cap[T(h_{\alpha+1}, f_{\alpha+1}, m')]\neq\emptyset$$
Without loss of generality $\alpha$ is even. Let $y$ be an element
of the first set and $y'$ an element of the second set.
By (iii) $y<^* f_{\alpha+1}$ and, by the remark after the construction
of $T(f,g,m)$, $y\not<^*f_\alpha$. Thus $y\in B$. Similarly $y'\in A$.
Hence the $\Dl{1}{2}(a)$--set $A$ is not $\MMM$--measurable,
a contradiction. This completes the proof.\qed

%%%%%%%%%%%%%%%%%%%%%%%%%%%   Section 7   %%%%%%%%%%%%%%%%%

\section{{\sc Sacks} Measurability} \label{Sackssec}
In this last section we will prove no new theorem, but apply well--known
results  to get an analogous characterization for
{\sc Sacks} measurability.
\begin{sthm}
The following are equivalent:
\begin{enumerate}
\item $\forall a\in\omom : \omom\cap{\bf L}[a]\neq\omom$
\item $\DiiS$
\item $\SiiS$
\end{enumerate}
\end{sthm}
{\bf Proof :} \\
For the direction (ii)$\Rightarrow$(i) we use the well-known construction
of a $\Dl{1}{2}(a)$--{\sc Bernstein} set in ${\bf L}[a]$. This leaves only
the direction (i)$\Rightarrow$(iii) to be proved. For this we need the 
theorem of {\sc Mansfield} and {\sc Solovay} ({\it cf.} 
\cite{S69} or \cite{Jech}, p. 533{\it sq.}):
\begin{sthm}[{\sc Mansfield--Solovay}]
If $A$ is a $\Sl{1}{2}(a)$--set of reals, then either it is in ${\bf L}[a]$ 
or it contains the branches through a perfect tree.
\end{sthm}
Now we take a $\Sl{1}{2}(a)$--set $A$ and show that it is $\SSS$--measurable.
As we have a real which is not constructible from $a$ we can even find
a perfect set of reals $P$ in $\omom\setminus{\bf L}[a]$.
{\sc Mansfield--Solovay} tells us that either $A$ contains a perfect set
or $A$ is in ${\bf L}[a]$ in which case $P$ lies completely in the
complement of $A$. Therefore we have $w\SiiS$ and with \ref{weakM}
even $\SiiS$.\qed

%%%%%%%%%%%%%%%%%%%%%%%%   Section 8   %%%%%%%%%%%%%%%%%%%%%%%

\section{Summary and Questions} \label{Summ}

We summarize our results and the older results of {\sc Solovay}
and \cite{JSh89} in the following table, where $\BBB$ denotes
Random Forcing, ${\sf Ran}(\frnM)$ the set of all random reals
over $\frnM$ and $\lambda$ the {\sc Lebesgue} measure:\\[5mm]
\begin{center}
\begin{tabular}{|c||c|c|}
\hline
\rule[-3mm]{0mm}{8,5mm}
Forcing & $\Db{1}{2}(\PPP)$ & $\Sb{1}{2}(\PPP)$\\ \hline\hline
\rule[-3mm]{0mm}{8,5mm}
$\BBB$ & $\forall a\in\omom : \RanL{a}\neq\emptyset$ & 
$\forall a\in\omom : \lambda(\RanL{a})=1$\\ \hline
\rule[-3mm]{0mm}{8,5mm}
$\CCC$ & $\forall a\in\omom : \CohL{a}\neq\emptyset$ & 
$\forall a\in\omom : \CohL{a}\in (c^1)$\\ \hline
\rule[-3mm]{0mm}{8,5mm}
$\DDD$ & $\forall a\in\omom : \HechL{a}\neq\emptyset$ & 
$\forall a\in\omom : \HechL{a}\in (d^1)$\\
\rule[-3mm]{0mm}{8,5mm}
& $\iff\SiiC$ & 
$\iff\forall a\in\omom : \aleph_1^{{\bf L}[a]} < \aleph_1$\\ \hline
\rule[-3mm]{0mm}{8,5mm}
$\LLL$ & \multicolumn{2}{c|}{$\forall a\in\omom :
\omom\cap{\bf L}[a]$ is $\sigma$--bounded}\\ \hline
\rule[-3mm]{0mm}{8,5mm}
$\MMM$ & \multicolumn{2}{c|}{$\forall a\in\omom :
\omom\cap{\bf L}[a]$ is not dominating}\\ \hline
\rule[-3mm]{0mm}{8,5mm}
$\SSS$ & \multicolumn{2}{c|}{$\forall a\in\omom :
\omom\cap{\bf L}[a]\neq\omom$}\\ \hline
\end{tabular}

\end{center}

By the  characterizations in the table and by well--known forcing arguments,
none of the arrows in Corollary \ref{Gamma} reverses (in ZFC)
for $\Gamma$ being either $\Db {1} {2}$ or $\Sb{1}{2}$.
However, for $\Gamma = \Sb{1}{2}$, the diagram gets
simpler because we then have $\Gamma(\CCC) \Longrightarrow
\Gamma (\LLL)$.

A few comments about full projective measurability in each of
our cases are in order. First, standard arguments show that
$\Sb{1}{n}(\PPP)$ holds, for all $n$ and all $\PPP$ considered
in this work, in {\sc Solovay}'s model which is
gotten by collapsing an inaccessible ({\it cf.}
\cite{S70} or \cite{Jech}, p. 537{\it sqq.}). Hence the consistency strength
of full projective measurability is at most an inaccessible.
In the {\sc Hechler} case, it is exactly an inaccessible by
\ref{HechSigma}.

Furthermore, $\Sb{1}{n}(\SSS)$ holds for all $n$ in
the model gotten by adding $\aleph_1$ {\sc Cohen} reals. To see this,
simply note that {\sc Cohen} forcing adds a perfect set of {\sc Cohen}
reals, and then use homogeneity of {\sc Cohen} forcing. 
Finally, $\Sb{1}{n}(\MMM)$ holds for all $n$ in {\sc Shelah}'s
model for the projective {\sc Baire} property ({\it cf.} \cite{Sh84}
or \cite{BJ95}, p. 495{\it sqq.}). This is true
by Corollary \ref{Gamma}. Hence in both cases the consistency
strength of full projective measurability is ZFC alone.
However, we do not know the answer to the following
\begin{sques}
Can one prove the consistency of ``all projective sets are
$\LLL$--measurable" on the basis of the consistency of ZFC
alone?
\end{sques}
Since {\sc Laver} forcing is closely related to {\sc Mathias} forcing,
this question has a flavour similar to the famous open
problem about the consistency strength of ``all projective
sets are completely {\sc Ramsey}".

\end{document}